# 基于协同多目标粒子群优化的交直流系统最优潮流

李亚辉,李扬,李国庆

(东北电力大学 电气工程学院,吉林省吉林市 132012)

**摘要:** 为了统一协调经济性和电压偏差,提出一种基于协同多目标粒子群优化(CMOPSO)的含 VSC-HVDC 交直流系统多目标最优潮流(MOPF)算法。首先,基于 VSC-HVDC 稳态模型,以最小化发电成本和电压偏差为目标,构建交直流系统的 MOPF 模型;然后,采用 CMOPSO 算法优化该模型,得到具有良好分布的 Pareto 最优解集;最后,在通过模糊 C 均值算法将所得解集进行聚类的基础上,采用灰关联投影法计算各决策方案的优属度,确定反映决策者不同偏好的最优折中解。应用于 IEEE 14 节点和 IEEE 118 节点系统的测试结果验证了所提方法的有效性。

**关键词:** 交直流系统;多目标最优潮流;柔性直流输电;协同多目标粒子群;多属性决策

## 0 引言

电力系统最优潮流(optimal power flow,OPF)是确保系统安全、经济运行的重要手段[1,2]。随着我国特高压交直流混联电网的快速发展、大规模清洁能源的接入,电力系统中电力电子器件广泛使用,系统运行特性更为复杂多变、日益逼近其稳定极限,传统的单一目标最优潮流已无法满足系统多目标协调优化运行的实际需求[3,4]。然而,多目标最优潮流(multi-objective OPF,MOPF)由于能够有效地求解包含重要性不同、甚至互相冲突的多个目标的问题,因而受到了国内外学者的关注[5,6]。同时,近年新兴的柔性直流输电(voltage source converter based high voltage direct current, VSC-HVDC)技术凭借其有功、无功功率可独立、快速控制、便于搭建多端直流(multi-terminal direct current,MTDC)网络等优点,在国内外大量的工程实践中获得了成功应用[7,8]。因此,求解含 VSC-HVDC 系统的 MOPF 成为亟待解决的重要问题。

国内外学者对含 VSC-HVDC 交直流系统的潮流优化问题已开展了大量研究。文献[9]采用两种改进内点法计算交直流系统的单目标最优潮流,并对结果进行了对比分析;文献[10]研究含多端柔直系统的最优潮流问题,对系统的经济性指标进行了优化;文献[11]基于二阶锥规划理论,采用内点法对单目标最优潮流问题进行讨论和求解;文献[12]以成本-效益指标为目标,对含 VSC-MTDC 交直流系统的最优潮流问题进行分析。然而,现有文献主要针对单目标 OPF 问题,无法满足智能电网建设中多元化的系统运行需求。电压偏差作为 MOPF 问题的常见指标,在传统交流系统中被广泛应用,但将该指标应用于交直流系统的研究却鲜有报道。而最近的研究表明,对于含 VSC-HVDC 交直流系统而言,维持足够的直流电压是运行过程中最重要的实际控制问题[13]。因此,本文在文献[14]的基础上从以下三个方面进行改进:在优化目标中,尝试加入电压偏差指标,从而减小电压的波动;在算法上,将多子群协同优化机制引入传统多目标进化算法中以加快其求解速度,缓解此类求解方法计算时间过长的问题,同时,考虑到运行人员对不同目标的偏好,采用聚类及决策分析的方法对解集进行评价;在系统中,增加直流网络 VSC 的数量,从两端直流网络扩展至多端直流网络,并对不同控制方式进行了对比。

本文提出一种基于协同多目标粒子群优化算法(cooperative multi-objective particle swarm optimization,CMOPSO)的含 VSC-HVDC 系统 MOPF 算法。首先,建立综合考虑发电成本和电压偏差的 MOPF 模型;然后,采用混合编码方式,通过 CMOPSO 得到该问题的 Pareto 最优解集;最后,在通过模糊 C 均值(fuzzy C-means,FCM)算法将

所得解集进行聚类的基础上,采用灰关联投影(grey relation projection,GRP)法评估各决策方案的优属度,从中选取最优折中解。

## 1 VSC-HVDC 稳态模型

### 1.1 稳态功率特性

含 VSC-HVDC 系统主要包含交流网络,直流网络和换流器三个部分。假设交直流系统中有 $N_{con}$ 个换流器,以第 $i$ 个换流器为例,则系统等效模型如下图所示。

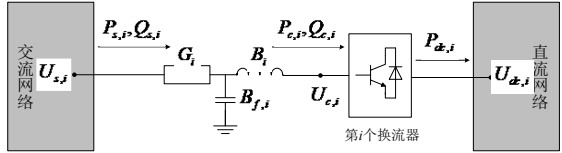

**图 1 交直流系统的简化等效模型图**
**Fig.1 Simplified equivalent model of AC/DC system**

图 1 中,换流变压器等效导纳可写为 $Y_i = G_i + jB_i$,其中低通滤波器 $B_{f,i}$ 的损耗可忽略,交流网络侧电压 $\dot{U}_{s,i}$、$U_{s,i} \angle \delta_{s,i}$,直流网络换流器输入电压 $\dot{U}_{c,i}$、$U_{c,i} \angle \delta_{c,i}$,则交流网络侧流出的有功功率 $P_{s,i}$ 及无功功率 $Q_{s,i}$ 分别为[15,16]:

$$P_{s,i} = U_{s,i}^2 G_i - U_{s,i}U_{c,i}\left[G_i \cos(\delta_{s,i} - \delta_{c,i}) + B_i \sin(\delta_{s,i} - \delta_{c,i})\right]$$
$$Q_{s,i} = -U_{s,i}^2 B_i - U_{s,i}U_{c,i}\left[G_i \sin(\delta_{s,i} - \delta_{c,i}) - B_i \cos(\delta_{s,i} - \delta_{c,i})\right]$$
(1)

由交流网络注入直流换流器的有功功率 $P_{c,i}$、无功功率 $Q_{c,i}$ 分别为:

$$P_{c,i} = -U_{c,i}^2 G_i + U_{s,i}U_{c,i}\left[G_i \cos(\delta_{s,i} - \delta_{c,i}) + B_i \sin(\delta_{s,i} - \delta_{c,i})\right]$$
$$Q_{c,i} = U_{c,i}^2 B_i - U_{s,i}U_{c,i}\left[G_i \sin(\delta_{s,i} - \delta_{c,i}) + B_i \cos(\delta_{s,i} - \delta_{c,i})\right]$$
(2)

注入直流换流器的有功功率 $P_{c,i}$ 与注入直流网络的有功功率 $P_{dc,i}$ 有如下关系:

$$P_{c,i} + P_{dc,i} + P_{con\_loss,i} = 0$$
$$s.t. \quad P_{dc,i} = U_{dc,i} \cdot I_{dc,i}$$
(3)

其中,$U_{dc,i}$、$I_{dc,i}$ 分别直流节点电压和电流,$P_{con\_loss,i}$ 为换流器损耗,表达式如下:

$$\begin{cases} P_{con\_loss,i} = a + b \cdot I_{c,i} + c \cdot I_{c,i}^2 \\ I_{c,i} = \dfrac{\sqrt{P_{c,i}^2 + Q_{c,i}^2}}{\sqrt{3}U_{c,i}} \end{cases}$$
(4)

式中,$I_{c,i}$ 为换流器电流;$a$、$b$ 及 $c$ 为各换流器相关系数,参考文献[17]各系数的取值如下: $a = 11.033 \times 10^{-3}$,$b = 3.464 \times 10^{-3}$,$c = 5.534 \times 10^{-3}$。

### 1.2 直流网络模型

对于两端直流系统,直流网络仅包含一条直流线路,模型相对简单。本节以三端直流网络为例,对直流网络模型进行分析,该直流网络的等效模型见附录 A 图 A1。

在对直流网络稳态模型分析时,线路中的功率仅与支路电阻以及直流节点的电压差有关,因此采用纯电阻电路对直流网络进行建模,注入节点 $i$ 的直流电流为:

$$I_{dc,i} = \sum_{j=1,j\neq i}^{N_{dc}} Y_{dc,ij}(U_{dc,i} - U_{dc,j})$$
(5)

式中,$Y_{dc,ij}$ 为节点 $i$ 和 $j$ 间的导纳,$N_{dc}$ 为直流网络节点数,该节点数与换流器个数 $N_{con}$ 相等。由上式可以看出,节点直流电压直接影响直流网络功率,因此为了维持功率的稳定,需要控制其直流电压。

### 1.3 换流器功率约束

对于含 VSC-HVDC 的交流系统,为保证系统的稳定运行,换流器需要运行在由电压 $U_{c,i}$ 和电流 $I_{c,i}$ 所确定的 P-Q 功率圆范围内,从而描述了换流器的有功和无功功率约束[17],第 $i$ 个换流器的 P-Q 功率圆见附录 A 图 A2,其表达式如下:

$$(r_i^{\min})^2 \leq (P_{s,i} - P_{0,i})^2 + (Q_{s,i} - Q_{0,i})^2 \leq (r_i^{\max})^2$$
(6)

上式表示换流器 $i$ 的功率圆,该圆以 $(P_{0,i}, Q_{0,i})$ 为圆心,以 $r_i$ 为半径,$r_i^{\max}$、$r_i^{\min}$ 分别为其上、下限。

由此可知,换流器输入电压、电流直接影响换流器功率,进而影响直流网络的传输功率及交直流系统的潮流分布[17],因此减小直流电压偏差对直流网络的运行十分重要。

### 1.4 稳态控制方式

VSC-HVDC 的基本控制方式可分为四种[10]: 1) 定 $U_{dc}$、定 $Q_s$ 控制;2) 定 $U_{dc}$、定 $U_s$ 控制;3) 定 $P_s$、定 $Q_s$ 控制;4) 定 $P_s$、定 $U_s$ 控制。

对于两端 VSC-HVDC 系统的控制方式,可采用 1)+3)、1)+4)、3)+2)、4)+2)。而当系统扩展为 VSC-MTDC 时,其控制方式可以派生出更多组合。但是在 VSC-MTDC 采用传统控制方式时,其中一端必须保持直流电压 $U_{dc}$ 不变,以维持功率平衡。随着技术的发展,VSC-MTDC 已产生电压下垂控制等新型的控制方式[18]。但以上控制方式均是为了防止直流电压 $U_{dc}$ 偏差过大,当系统直流偏差过大或当系统发生故障时,会造成交直流系统中潮流较大的波动[19]。因此控制系统的直流电压偏差,能够有效保证系统的稳定性[13,18]。

## 2 问题描述

### 2.1 目标函数

为了统一协调系统运行的经济性及电压偏差,采用发电成本最小及电压偏差最小作为优化目标。

### 2.1.1 系统发电成本

将发电成本 $F$ 最小作为目标可以有效地优化系统有功功率，其表达式如下：

$$F = \sum_{i=1}^{N_G}(\alpha_i P_{G,i}^2 + \beta_i P_{G,i} + \gamma_i) \quad (7)$$

式中，$N_G$ 为发电机台数，$P_{G,i}$ 为第 $i$ 台发电机的有功出力，$\alpha_i$、$\beta_i$ 和 $\gamma_i$ 分别为发电机 $i$ 的各项发电成本系数。

### 2.1.2 电压偏差指标

本文将电压偏差指标 $V_{de}$ 最小作为另一个优化目标：

$$V_{de} = \sum_{i=1}^{N}(U_i - U_{ref,i})^2 + \sum_{j=1}^{N_{dc}}(U_{dc,j} - U_{ref,dc,j})^2 \quad (8)$$

式中，$U_i$ 为交流网络节点 $i$ 的电压，其给定参考值为 $U_{ref,i}$；$U_{dc,i}$ 为直流网络节点 $j$ 的电压，其给定参考值为 $U_{ref,dc,i}$；$N$ 为交流网络节点数。

在对交直流系统潮流进行优化的过程中，直流电压 $U_{dc}$ 作为控制变量，能够通过调整电压 $U_c$ 使换流器工作在 $P$-$Q$ 功率圆范围内，而且能够保证直流功率的稳定；同时，考虑到 VSC 的控制方式中，需要对直流电压 $U_{dc}$ 进行控制，以防止电压偏差过大[18]。基于以上两方面，在电压偏差指标 $V_{de}$ 的表达式中，加入直流电压 $U_{dc}$ 的偏差，在一定程度上保证了系统稳定运行。

## 2.2 约束条件

本节对含 VSC-HVDC 交直流系统中相关约束条件进行介绍。

交流网络中等式约束为节点功率平衡方程；不等式约束包括以下变量的限制：发电机有功、无功出力，节点电压，可调变压器分接头变比，无功补偿装置投切容量，线路传输容量。因为篇幅所限，交流网络约束的具体公式在此不再赘述，详见参考文献[14]。

在直流网络和换流站中，相关等式约束已在1.1 节中进行介绍；不等式约束除 1.2 节中介绍的换流器功率约束外，还包括以下变量的约束：

$$\begin{cases} P_{s,i}^{\min} \leq P_{s,i} \leq P_{s,i}^{\max}, & i=1,\ldots N_{dc} \\ Q_{s,i}^{\min} \leq Q_{s,i} \leq Q_{s,i}^{\max}, & i=1,\ldots N_{dc} \\ U_{dc,i}^{\min} \leq U_{dc,i} \leq U_{dc,i}^{\max}, & i=1,\ldots N_{dc} \\ I_{dc,i}^{\min} \leq I_{dc,i} \leq I_{dc,i}^{\max}, & i=1,\ldots N_{dc} \\ -I_{dc,ij}^{\max} \leq I_{dc,ij} \leq I_{dc,ij}^{\max}, & i,j=1,\ldots N_{dc}\ i \neq j \end{cases} \quad (9)$$

式中，上标 max 和 min 分别代表对应物理量的上、下限，$I_{dc,ij}$ 为直流节点 $i$、$j$ 之间流过的电流。

## 3 模型求解

### 3.1 求解框架

MOPF 既是优化问题，也是决策问题[16]，因此求解过程分为多目标优化和决策支持两部分。在多目标优化部分，采用 CMOPSO 算法对求解所建 MOPF 模型，得到 Pareto 最优解集；在决策支持部分，对所得解集使用 FCM 聚类及 GRP 相结合的方式，通过进一步决策分析给出最优折中解。整体求解框架如图 2 所示。

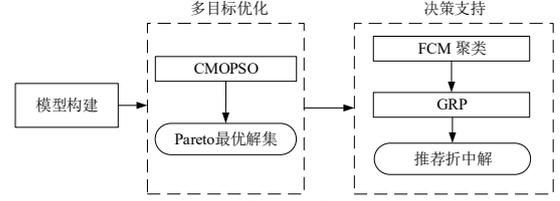

图 2 整体求解框架
Fig.2 Overall framework of the proposed approach

### 3.2 多目标优化

多目标粒子群优化算法是求解多目标优化问题的有力工具，具有概念简单、易于实现、收敛速度快等优点[20]。由于 OPF 问题对求解速度要求较高，为此提出融合了多子群协同优化机制的 CMOPSO 算法，进一步提高计算速度，该算法的主要步骤为：

1）输入系统参数、算法参数和控制变量参数，各类参数分别包含以下数据：

①计算潮流所需的系统数据，对于直流网络，需指定 VSC 的数目及其控制方式；

②种群规模 $S_{pop}$、子种群规模 $S_{sub}$、最大迭代次数 $I_{\max}$，交换间隔次数 $I_t$，迭代次数 $N_{ite}=0$，目标函数及待优化变量个数；

③待优化变量的范围，其中，离散变量需要给定步长。

2）初始化粒子位置和速度，并用混合编码方式对决策变量进行编码。编码时，连续变量采用实数编码，包含交流网络侧流出的有功功率 $P_s$、无功功率 $Q_s$ 及交流侧电压 $U_s$ 等直流量及发电机有功出力 $P_G$、机端电压 $U_G$ 等交流量；离散变量采用整数编码，包括可调变压器分接头变比 $T$ 和无功补偿电容器投切容量 $Q_C$。

3）依据粒子位置和速度，对交直流系统数据进行更新，进而采用交替迭代法[17]，计算交直流系统潮流；由交流系统潮流结果，可得各粒子的目标函数向量 $\boldsymbol{f}=(F,V_{de})$；通过目标函数值的比较，能够得到其中的非支配解，并建立外部档案进行存储。在此过程中，各粒子需要记录自身的个体历史最优位置 $pbest$，同时找出整个种群及各个子种群的全局最优位置 $gbest$ 及 $sub\text{-}gbest$。

4）按照算法自身寻优机制，分别对各子种群进行优化，确定新的粒子位置及速度，进而形成新的子种群。各子种群之间采用协同优化方式，运行过程中同时对多个子群进行处理。

5）计算新形成子种群的目标函数向量 $\boldsymbol{f}$，更新外部档案及各粒子的历史最优位置 $pbest$，同时找出全局最优位置 $gbest$ 及 $sub\text{-}gbest$。

6) 若达到交换间隔次数 $I_t$，则交换相邻子种群的全局最优位置 sub-gbest；否则，继续执行程序。

7) 终止条件判断：若 $N_{ite} = I_{max}$，则输出 Pareto 最优解集；否则，$N_{ite} = N_{ite} + 1$，更新数据后转至步骤 4)继续执行程序。

### 3.3 决策支持

优化所得 Pareto 最优解规模大，且决策向量包含不同信息，本文提出 FCM 聚类与 GRP 相结合的辅助决策方法，便于运行人员从中选取折中解。

#### 3.3.1 模糊 C 均值聚类

FCM 是一种非监督聚类算法，其数学模型为：

$$\min J_n(S,C) = \sum_{ii=1}^{N_p}\sum_{jj=1}^{N_{clu}} \eta_{ii,jj}^n \|s_{ii} - c_{jj}\|^2 \quad (10)$$

$$s.t. \sum_{jj=1}^{N_{clu}} \eta_{ii,jj} = 1$$

其中，$J$ 为聚类损失函数，$S = \{s_1, s_2, \cdots \cdots\}$ 为 Pareto 最优解的集合，$N_p$ 为集合中解的个数，$C = \{c_1, c_2, \cdots \cdots\}$ 为各聚类中心，$N_{clu}$ 为预先设定的聚类数目，$\eta_{ii,jj}(\eta_{ii,jj} \in [0,1])$ 为第 $ii$ 个样本对第 $jj$ 类的隶属度，$n(n \in [1,\infty])$ 为控制聚类模糊程度常数。

采用 FCM 聚类后，可使同类解之间相似度最大，而不同类解之间相似度最小。由于本文中考虑两个目标，为反映运行人员对目标的不同偏好，因此将聚类数目设为 2。

#### 3.2.2 灰关联投影法

通过将灰色系统理论和矢量投影原理结合，获得适合处理实际系统中的灰色多属性决策问题的灰关联投影法[21,22]。由指标特性可知，文中所用两个目标均为效益型指标，其中第 $l$ 个方案在理想方案上的投影值 $V_l^{+(-)}$ 为：

$$V_l^{+(-)} = \sum_{k=1}^{t} \gamma_{lk}^{+(-)} \frac{w_{G,k}^2}{\sqrt{\sum_{k=1}^{t}(w_k)^2}} \quad (11)$$

式中，上标"+"表示理想方案，上标"−"表示负理想方案；$t$ 为指标总数；$\gamma_{lk}^{+(-)}$ 为第 $l$ 个方案第 $k$ 个指标所对应与理想（负理想）方案的灰关联系数；$w_{G,k}$ 为方案各指标权重，文中将三个目标所对应权重设为相同数值，运行人员可根据实际工况或个人偏好进行调整。同时，定义优属度 $d$ 为：

$$d_l = \frac{(V_0 - V_l^-)^2}{(V_0 - V_l^-)^2 + (V_0 - V_l^+)^2}, \quad 0 \leq d_l \leq 1 \quad (12)$$

式中，$V_0$ 为 $V_l$ 在 $\gamma$ 取 1 时的值。由上式可知，为使所选方案更加接近理想方案，同时远离负理想方案，因此本文所提方法选择优属度最大的方案作为最优折中解。

## 4 算例分析

### 4.1 IEEE 14 节点系统

为验证所提算法的有效性，分别以含两端、三端直流网络的 IEEE 14 节点系统为例进行分析。

#### 4.1.1 含两端直流网络的 IEEE 14 节点系统

（1）算例介绍

该系统包括 2 台可控发电机，18 条支路，11 个负荷，4-5 支路修改为直流支路[9]，见附录 A 图 A3。设定各节点电压范围为[0.94, 1.06]；可调变压器分接头变比为[0.9, 1.1]，其步长为 0.0125；节点 9 的无功补偿电容器投切容量 $Q_{C,1}$ 范围为[0, 0.5]，其步长为 0.01；直流网络有功、无功功率范围为[-1.0, 1.0]。

直流节点参数见附录 A 表 A1，$B_{vsc}$ 表示 VSC 所连交流母线的编号；其中初始 $P_s$、$Q_s$ 及 $U_{dc}$ 设定依据未改造的交流系统 OPF 结果。在计算过程中，由于直流网络 VSC 控制方式的不同，其中的部分变量将作为控制变量而发生数值上的变化；同时，直流支路参数与修改前交流直流参数相同。

本节中 $VSC_1$ 和 $VSC_2$ 分别采用 1.4 节中介绍的控制方式 1) 和 3)。需要说明的是，本节仅以一种控制方式的组合为例进行分析，但该方法同样适用于其他控制方式。

（2）优化结果

为了合理评价所提 CMOPSO 的优化性能，将其与非支配遗传算法（non-dominated sorting genetic algorithm-II，NSGA-II）在优化效果和计算速度两方面进行了对比。其中，NSGA-II 由于具有求解速度快、收敛性好等优点，已被广泛用于求解工程领域中的各种复杂优化问题[14,16]。经过 50 次迭代，两种算法所得 Pareto 最优解集在目标空间上的分布见附录 A 图 A4，由图可知，虽然 NSGA-II 也可获得大体完整、均匀分布的 Pareto 最优解集，但所提算法优化效果明显优于 NSGA-II，较之更加接近于 Pareto 前沿。所得极端解见附录 A 表 A2，由表可知，针对各目标函数，所提算法所得极端解均优于 NSGA-II，进一步印证了其优越性。

为了分析两种算法的求解速度，将二者各自独

立重复运行 20 次的平均时间进行对比，结果见附录 A 表 A3。所用电脑配置为 CPU 主频 3.2 GHz，内存 4 GB。由表可知，所提算法的求解速度优于 NSGA-II。因此，CMOPSO 的优化效果和求解速度均优于常用的 NSGA-II 算法。

所得 Pareto 解集的 FCM 聚类结果，见附录 A 图 A5。其中，红色及绿色分别代表决策者偏好系统发电成本 $F$ 及电压偏差指标 $V_{de}$。

采用 GRP 进行方案评估时，将聚类后各类解作为决策单元，计算其优属度，取各类解中优属度最大的解作为最优折中解，结果见表 1。

表 1 IEEE 14 节点两端算例最优折中解
Table 1 Compromise solutions of 2-terminal DC network

| 最优折中解 | $F$ /($/h) | $V_{de}$ /p.u. | 优属度 |
|---|---|---|---|
| 折中解 1 | 8172.08 | 0.0095 | 0.7770 |
| 折中解 2 | 8187.80 | 0.0051 | 0.6403 |

由此可知，所提算法通过有效结合 CMOPSO 和 GRP，能够求解该算例的 MOPF 问题，并筛选出最优折中解。

（3）潮流优化前后结果对比

以表 1 折中解 1 为例，系统潮流优化前后的结果对比见附录 A 表 A4~6，由表可知，控制变量均在设定的优化范围内进行调整，通过潮流优化和直流系统的调节作用，使系统发电成本 $F$ 及电压偏差指标 $V_{de}$ 得以改善，因此优化后的潮流分布更为合理。至此，通过含两端直流网络的 IEEE 14 节点系统的分析，初步验证了所提算法的有效性。

（4）直流电压偏差对优化结果的影响

为了验证将直流电压偏差加入目标函数 $V_{de}$ 中的必要性，对分别采用目标函数 $V_{de}$ 和未加入直流电压偏差的目标函数 $V_{de,ac}$ 进行多目标潮流优化的结果进行对比分析，其中，目标函数 $V_{de,ac}$ 的表达式为：

$$V_{de,ac} = \sum_{i=1}^{N}(U_i - U_{ref,i})^2 \qquad (13)$$

算例中保持 VSC$_1$、VSC$_2$ 控制方式不变，Case 1 的优化目标为 $F$ 和 $V_{de}$，Case 2 的优化目标为 $F$ 和 $V_{de,ac}$。优化结束后，选取各工况解集中优属度最大的解作为参考解，并根据其运行点分别计算 Case 1 和 Case 2 的 $V_{de,ac}$ 和 $V_{de}$，算例的对比结果见附录 A 表 A7。由表中对比结果可以看出，虽然 Case 2 采用 $V_{de,ac}$ 作为目标函数进行优化，但所得目标函数值更大，但是 Case 2 在优化过程中，有可能会导致直流电压偏差过大，接近直流电压极限，不能对电压偏差进行有效的控制，影响系统稳定运行。因此可以得出结论，将直流电压偏差加入目标函数中是有必要的，可以保证直流电压稳定在设定范围内，维持系统中潮流的稳定。

4.1.2 含三端直流网络的 IEEE 14 节点系统

为了验证所提算法对含多端直流网络交直流系统的有效性，以 IEEE 14 节点三端交直流系统为例[16]，见附录 A 图 A6。直流节点参数见附录 A 表 A8。本节中 VSC$_3$ 采用 1.4 节中介绍的控制方式 1)，VSC$_1$ 和 VSC$_2$ 采用控制方式 3)。对于其他控制方式的组合，本文已进行检验，同样可以采用该方法进行求解。

所得 Pareto 最优解集经过 FCM 聚类后，在目标函数空间上的分布见附录 A 图 A7。采用 GRP 法对进行评估，所得最优折中解如表 2 所示。

表 2 IEEE 14 节点三端算例最优折中解
Table 2 Recommended solutions of 3-terminal DC

| 最优折中解 | $F$ /($/h) | $V_{de}$ /p.u. | 优属度 |
|---|---|---|---|
| 折中解 1 | 8143.64 | 0.0067 | 0.7805 |
| 折中解 2 | 8160.62 | 0.0040 | 0.7082 |

因此，本文所提算法对于含多端直流网络的交直流系统同样适用。

4.1.3 控制方式及 VSC 数量对潮流优化的影响

为了分析控制方式及 VSC 数量对潮流优化的影响，采用以下三种典型工况进行对比测试：

Case 0—未改造的 IEEE 14 节点系统，其中，直流电压偏差为 0；

Case 1—IEEE 14 节点两端交直流系统，VSC$_1$ 和 VSC$_2$ 分别采用 1.4 节中介绍的控制方式 1）和 3）；

Case 2—IEEE 14 节点两端交直流系统，对调 Case1 中的 VSC$_1$、VSC$_2$ 控制方式；

Case 3—IEEE 14 节点三端交直流系统，VSC$_1$ 采用 1.4 节中介绍的控制方式 1)，VSC$_2$ 和 VSC$_3$ 采用控制方式 3)；

Case 4—IEEE 14 节点三端交直流系统，VSC$_2$ 采用 1.4 节中介绍的控制方式 1)，VSC$_1$ 和 VSC$_3$ 采用控制方式 3)；

Case 5—IEEE 14 节点三端交直流系统，VSC$_3$ 采用 1.4 节中介绍的控制方式 1)，VSC$_1$ 和 VSC$_2$ 采用控制方式 3)；

Case 6—IEEE 14 节点三端交直流系统，VSC$_1$、VSC$_2$ 及 VSC$_3$ 的控制方式均为下垂控制方式，其下垂斜率均设置为 0.005[17]。需要说明的是，Case 6 中可以将下垂斜率作为控制变量进行优化，但本文为了尽可能的保证控制变量的一致性，没有对下垂斜率进行优化。

采用 CMOPSO 算法优化后，以 Case 0~6 的解集中优属度最大解为例，优化后系统部分参数

见附录 A 表 A9，其中，$IMP_F$、$IMP_V$ 分别表示相对交流系统，发电成本 $F$ 和电压偏差指标 $V_{de}$ 的改善百分比。由表中数据，可做如下分析：1）通过 Case 0 和其他工况的对比可知，系统中加入 VSC 后两个目标函数的值变小，说明相比于纯交流系统，含 VSC-HVDC 的交直流系统具有较强的调节能力；2）通过对含有不同 VSC 数量的系统进行对比，可以发现，优化后多端直流网络的发电成本及电压偏差指标均优与两端直流网络。这说明增加系统中 VSC 的数量可以使系统潮流更为合理，也说明直流系统的调节能力随 VSC 数量增加而增强；3）对于具有相同拓扑结构的交直流系统而言，VSC 控制方式的选取对潮流优化结果有一定的影响：Case 1 的发电成本 $F$ 优于 Case 2，但 Case 2 的电压偏差指标 $V_{de}$ 优于 Case 1，而根据 Case 3~6 中目标函数值，各算例并非比其他算例更具优势，且目标函数在数值上差别不大，说明控制策略的改变对优化结果有一定影响，但不一定改善优化目标。

## 4.2 IEEE 118 节点系统

### 4.2.1 含三端直流网络的 IEEE 118 节点系统

为考察所提算法对较大规模系统的适用性，采用经修改的 IEEE 118 节点系统进行仿真验证。该系统包括有 14 台可控有功出力发电机；并分别将 $VSC_1$，$VSC_2$ 和 $VSC_3$ 连接系统中 103、105 和 104 节点，阻抗参数同附录 A 表 A8，$P_s$、$Q_s$ 及 $U_{dc}$ 参考未改造系统的 OPF 结果取得；系统中三端直流网络控制方式同 4.1.2 节。

经过 100 次迭代，所得 Pareto 最优解集聚类后在目标函数空间上的分布见附录 A 图 A8。

可知，对于 IEEE 118 节点系统，所提算法仍可得到分布较为完整且均匀的 Pareto 最优解集。经 GRP 方案评估，所得最优折中解见表 3。

表 3 IEEE 118 节点算例最优折中解
Table 3 Recommended solutions of IEEE 118-bus

| 最优折中解 | $F$ /($/h) | $V_{de}$/p.u. | 优属度 |
|---|---|---|---|
| 折中解 1 | 130906 | 0.0308 | 0.7044 |
| 折中解 2 | 131097 | 0.0223 | 0.6802 |

由此可知，所提算法能够获得 118 节点算例的最优折中解，因此同样适用于 IEEE 118 节点系统，从而验证了其对较复杂系统的适用性。

### 4.2.2 VSC 数量对潮流优化的影响

为了进一步分析较大规模系统中 VSC 数量对系统潮流优化问题的影响，采用以下三个算例进行对比：Case 0—未改造的 IEEE 118 节点系统；Case 1—IEEE 118 节点两端系统，其中，$VSC_1$ 和 $VSC_2$ 分别 103、104 节点相连，$VSC_1$ 和 $VSC_2$ 分别采用 1.4 节中介绍的控制方式 3）和 1）；Case 2—IEEE 118 节点三端系统，控制方式同 4.2.1 节；Case 3—IEEE 118 节点三端系统，控制策略修改为下垂控制方式，其下垂斜率均设置为 0.005[17]。采用 CMOPSO 算法优化后，以 Case 0~3 的解集中优属度最大解为例，优化后系统部分参数见附录 A 表 A10。限于文章篇幅，表中仅给出直流节点相关参数。由表中数据，可知，相对于交流网络，系统中加入 VSC 后，优化效果更好；相比于两端直流网络，多端直流网络的目标函数在优化后更具优势，可以得出与 4.1.3 节中相似的结论。但是，由于直流网络相对于交流网络而言规模较小，优化目标值没有十分显著的改善，说明交流系统的规模对直流网络的调控能力有一定的影响，当系统规模较大时，难以发挥 VSC 的调节能力。

## 5 结语

提出了一种基于 CMOPSO 的含 VSC-HVDC 系统 MOPF 算法，并以 IEEE 14 节点、118 节点系统为例进行了仿真研究，结论如下：

1）构建了计及发电成本和电压偏差的交直流系统 MOPF 模型，统一协调系统运行的经济性及电压偏差，适应了系统多目标协调优化运行的实际需求。

2）所提基于 CMOPSO 的模型求解方法，在求得 Pareto 最优解集的同时，可进一步从中筛选出最优折中解；且 CMOPSO 的优化性能优于 NSGA-II。

3）选用不同的控制方式及 VSC 数量对潮流优化结果具有一定影响。直流系统调节能力随着 VSC 数量增多而增强，潮流分布结果也随之更优。

此外，引入并行计算以提高求解速度、计及安全约束、VSC 控制方式优选及优化布点等问题是本文下一步工作要考虑的重点。

李亚辉(1993—)，男，硕士研究生，主要研究方向：柔性直流输电技术、电力系统经济运行等。E-mail：liyh1993@gmail.com

李扬(1980—)，男，通信作者，博士，副教授，主要研究方向：电力系统安全稳定分析与控制。E-mail：liyang@neppu.edu.cn

李国庆(1963—)，男，教授，博士生导师，主要研究方向：电力系统安全稳定分析与控制。E-mail：lgq@neepu.edu.cn

**Optimal Power Flow for AC/DC System Based on Cooperative Multi-objective Particle Swarm Optimization**

LI Yahui, LI Yang, LI Guoqing

(School of Electrical Engineering, Northeast Electric Power University, Jilin 132012, China）

**Abstract:** In order to unifiedly coordinate economy and voltage deviations, a novel multi-objective optimal power flow (MOPF) algorithm is proposed for an AC/DC system with VSC-HVDC based on cooperative multi-objective particle swarm optimization (CMOPSO). In order to minimize power generation costs and voltage deviations, the MOPF model of the AC/DC system is firstly built based on the VSC-HVDC steady-state model. Then, the CMOPSO is adopted for solving the MOPF model to find well-distributed Pareto-optimal solutions. Next, the solutions are divided into different groups via the fuzzy C-means algorithm, and finally the best compromise solutions reflecting decision-makers' different preferences are identified from each group by comparing the priority memberships which are calculated by using the grey relation projection method. The validity of the proposed approach is verified by using the modified IEEE 14-bus and 118- bus systems.

This work is supported by National Key R&D Program of China (No. 2017YFB0902401) and National Natural Science Foundation of China (No. 51677023).

**Key words:** AC/DC system; multi-objective optimal power flow; VSC-HVDC; cooperative multi-objective particle swarm optimization; multiple attribute decision making

## 附录 A

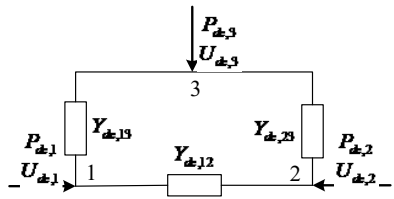

图 A1　三端直流网络模型图
Fig.A1　Model of three bus DC network

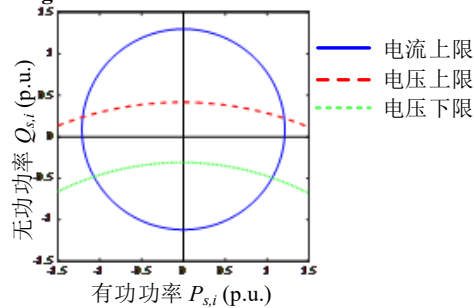

图 A2　第 *i* 个换流器 *P-Q* 功率圆
Fig.A2　*P-Q* capability chart of converter station *i*

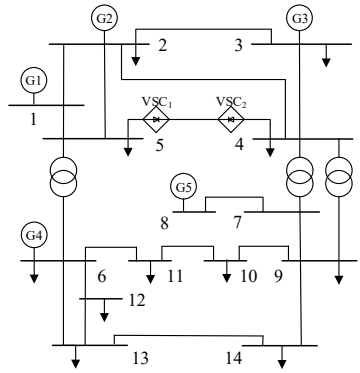

图 A3 修改的 14 节点两端交直流系统
Fig.A3 Modified 14-bus system with 2-terminal DC network

表 A1 两端直流系统节点参数
Table A1 Parameters of buses for 2-terminal DC network

| $B_{vsc}$ | $R$ /p.u. | $X$ /p.u. | $P_s$ /p.u. | $Q_s$ /p.u. | $U_{dc}$ /p.u. |
|---|---|---|---|---|---|
| 4 | 0.0015 | 0.1121 | -0.492 | 0.116 | 1.000 |
| 5 | 0.0015 | 0.1121 | 0.495 | -0.105 | 1.000 |

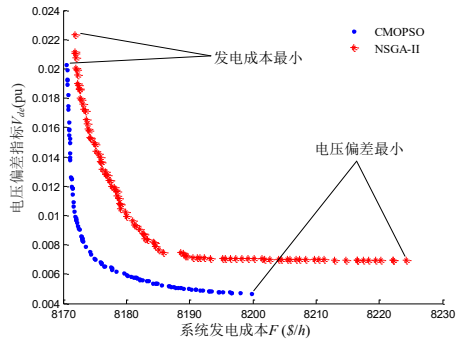

图 A4 IEEE 14 节点两端算例 Pareto 最优解的分布
Fig.A4 Distribution of solutions of IEEE 14-bus system with 2-terminal DC network

表 A2 IEEE 14 节点两端算例极端解
Table A2 Extreme solutions of 14-bus with 2-terminal DC

| 极端解 | 优化目标 | $F$ /($/h) | $V_{de}$ /p.u. |
|---|---|---|---|
| NSGA-II | 发电成本最小 | 8171.79 | 0.0223 |
|  | 电压偏差最小 | 8224.23 | 0.0069 |
| CMOPSO | 发电成本最小 | 8170.53 | 0.0203 |
|  | 电压偏差最小 | 8199.80 | 0.0047 |

表 A3 两种算法的求解时间
Table A3 Solution times of CMOPSO and NSGA-II

| 评价指标 | NSGA-II | CMOPSO |
|---|---|---|
| 平均时间(s) | 170.53 | 127.81 |

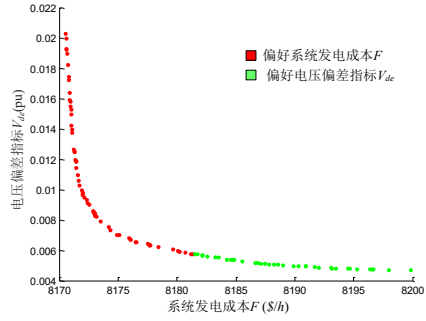

图 A5　IEEE 14 节点两端算例 Pareto 最优解分布
Fig.A5　Distribution of solutions of IEEE 14-bus system with 2-terminal DC network

表 A4　优化前后直流网络变量对比
Table A4　Comparison of DC network variables

| 直流网络 | 优化前 | | | 优化后 | | |
|---|---|---|---|---|---|---|
| | $P_s$ /p.u. | $Q_s$ /p.u. | $U_{dc}$ /p.u. | $P_s$ /p.u. | $Q_s$ /p.u. | $U_{dc}$ /p.u. |
| VSC$_1$ | 0.495 | -0.105 | 1.000 | 0.506 | -0.095 | 0.995 |
| VSC$_2$ | -0.492 | 0.116 | 1.000 | -0.504 | 0.104 | 1.049 |

表 A5　优化前后发电机变量对比
Table A5　Comparison of generator variables

| 发电机 | 优化前 | | | 优化后 | | |
|---|---|---|---|---|---|---|
| | $P_s$ /p.u. | $Q_s$ /p.u. | $U_{dc}$ /p.u. | $P_s$ /p.u. | $Q_s$ /p.u. | $U_{dc}$ /p.u. |
| G1 | 2.324 | -0.165 | 1.060 | 1.952 | -0.097 | 1.060 |
| G2 | 0.400 | 0.436 | 1.045 | 0.369 | 0.180 | 1.044 |
| G3 | 0 | 0.251 | 1.010 | 0.299 | 0.203 | 1.021 |
| G4 | 0 | 0.127 | 1.070 | 0.001 | 0.363 | 1.034 |
| G5 | 0 | 0.176 | 1.090 | 0.085 | 0.035 | 1.059 |

表 A6　优化前后目标函数对比
Table A6　Comparison of objective functions

| 优化状态 | $F$ /($/h) | $V_{de}$ /p.u. |
|---|---|---|
| 优化前 | 8287.68 | 0.0232 |
| 优化后 | 8172.08 | 0.0095 |

表 A7　不同目标函数算例优化结果对比
Table A7　Comparison of results with different objectives

| 变量 | Case 1 | Case 2 |
|---|---|---|
| $U_1$ /p.u. | 1.060 | 1.060 |
| $U_2$ /p.u. | 1.044 | 1.042 |
| $U_3$ /p.u. | 1.012 | 1.016 |
| $U_4$ /p.u. | 1.027 | 1.027 |
| $U_5$ /p.u. | 1.024 | 1.024 |
| $U_6$ /p.u. | 1.028 | 1.025 |
| $U_7$ /p.u. | 1.041 | 1.040 |
| $U_8$ /p.u. | 1.054 | 1.053 |
| $U_9$ /p.u. | 1.021 | 1.020 |
| $U_{10}$ /p.u. | 1.015 | 1.013 |
| $U_{11}$ /p.u. | 1.018 | 1.015 |
| $U_{12}$ /p.u. | 1.013 | 1.010 |
| $U_{13}$ /p.u. | 1.009 | 1.006 |
| $U_{14}$ /p.u. | 0.997 | 0.995 |
| $U_{dc,1}$ /p.u. | 1.000 | 1.020 |
| $U_{dc,2}$ /p.u. | 0.997 | 1.018 |
| $F$ /($/h) | 8174.34 | 8174.05 |
| $V_{de}$ /p.u. | 0.007338 | 0.007962 |
| $V_{de,ac}$ /p.u. | 0.007327 | 0.007241 |

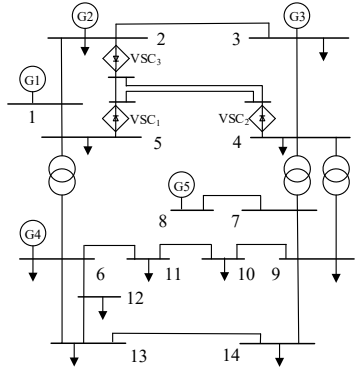

图 A6　修改的 IEEE 14 节点三端交直流系统
Fig.A6　Modified IEEE 14-bus with 3-terminal DC network

表 A8　三端直流系统节点参数
Table A8　Bus parameters of 3-terminal DC network

| $B_{VSC}$ | $R$ /p.u. | $X$ /p.u. | $P_s$ /p.u. | $Q_s$ /p.u. | $U_{dc}$ /p.u. |
|---|---|---|---|---|---|
| 2 | 0.0015 | 0.150 | 0.877 | 0.001 | 1.000 |
| 4 | 0.0015 | 0.150 | -0.983 | 0.124 | 1.000 |
| 5 | 0.0015 | 0.150 | 0.118 | -0.135 | 1.000 |

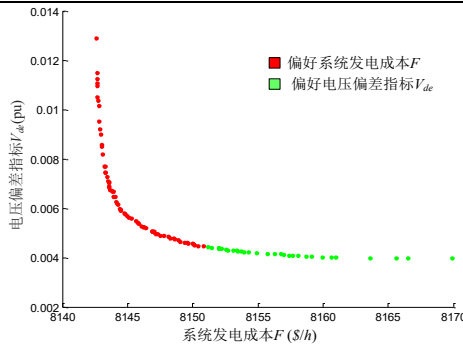

图 A7　14 节点三端算例 Pareto 最优解分布
Fig.A7　Distribution of solutions of 3-terminal DC network

表 A9　IEEE 14 节点各工况优化结果对比
Table A9　Comparison of optimal results of 14-bus system

| 变量 | 最大值 | 最小值 | Case 0 | Case 1 | Case 2 | Case 3 | Case 4 | Case 5 | Case 6 |
|---|---|---|---|---|---|---|---|---|---|
| $P_{G,1}$ /p.u. | 3.32 | 0 | 1.949 | 1.955 | 1.963 | 1.961 | 1.961 | 1.962 | 1.962 |
| $P_{G,2}$ /p.u. | 1.40 | 0 | 0.366 | 0.37 | 0.37 | 0.38 | 0.373 | 0.374 | 0.372 |
| $P_{G,3}$ /p.u. | 0.30 | 0 | 0.300 | 0.300 | 0.299 | 0.300 | 0.300 | 0.299 | 0.300 |
| $P_{G,4}$ /p.u. | 0.10 | 0 | 0.000 | 0.001 | 0.002 | 0.000 | 0.000 | 0.001 | 0.002 |
| $P_{G,5}$ /p.u. | 0.10 | 0 | 0.079 | 0.082 | 0.064 | 0.059 | 0.071 | 0.065 | 0.064 |
| $U_{G,1}$ /p.u. | 1.06 | 0.94 | 1.060 | 1.060 | 1.060 | 1.060 | 1.060 | 1.060 | 1.060 |
| $U_{G,2}$ /p.u. | 1.06 | 0.94 | 1.040 | 1.044 | 1.040 | 1.042 | 1.043 | 1.041 | 1.042 |
| $U_{G,3}$ /p.u. | 1.06 | 0.94 | 1.010 | 1.012 | 1.012 | 1.012 | 1.012 | 1.012 | 1.012 |
| $U_{G,4}$ /p.u. | 1.06 | 0.94 | 1.035 | 1.028 | 1.021 | 1.027 | 1.025 | 1.026 | 1.025 |
| $U_{G,5}$ /p.u. | 1.06 | 0.94 | 1.060 | 1.054 | 1.060 | 1.060 | 1.060 | 1.060 | 1.060 |
| $P_{s,1}$ /p.u. | 1.00 | -1.00 | - | 0.506 | 0.550 | 0.830 | 0.785 | 0.863 | -0.847 |

| | | | | | | | | |
|---|---|---|---|---|---|---|---|---|
| $Q_{s,1}$ /p.u. | 1.00 | -1.00 | - | -0.094 | -0.095 | 0.010 | 0.016 | 0.021 | 0.077 |
| $U_{dc,1}$ /p.u. | 1.06 | 0.94 | - | 1.001 | 1.000 | 0.989 | 1.000 | 0.997 | 0.985 |
| $P_{s,2}$ /p.u. | 1.00 | -1.00 | - | -0.505 | -0.523 | -0.949 | -0.891 | -0.983 | -0.986 |
| $Q_{s,2}$ /p.u. | 1.00 | -1.00 | - | 0.104 | 0.104 | 0.126 | 0.097 | 0.107 | 0.082 |
| $U_{dc,2}$ /p.u. | 1.06 | 0.94 | - | 0.997 | 0.997 | 1.001 | 0.986 | 1.011 | 0.999 |
| $P_{s,3}$ /p.u. | 1.00 | -1.00 | - | - | - | 0.131 | 0.118 | 0.132 | 0.126 |
| $Q_{s,3}$ /p.u. | 1.00 | -1.00 | - | - | - | -0.117 | -0.049 | -0.044 | -0.044 |
| $U_{dc,3}$ /p.u. | 1.06 | 0.94 | - | - | - | 0.998 | 0.998 | 1.000 | 0.996 |
| $T_1$ /p.u. | 1.1 | 0.9 | 0.9780 | 0.9655 | 0.9655 | 1.0405 | 1.0405 | 1.0405 | 1.0530 |
| $T_2$ /p.u. | 1.1 | 0.9 | 0.9690 | 1.0690 | 1.0565 | 0.9815 | 0.9690 | 0.9690 | 0.9315 |
| $T_3$ /p.u. | 1.1 | 0.9 | 0.9320 | 1.0320 | 1.0320 | 0.9820 | 0.9320 | 0.9320 | 0.9570 |
| $Q_{C,1}$ /p.u. | 0.5 | 0 | 0.21 | 0.19 | 0.19 | 0.20 | 0.17 | 0.19 | 0.22 |
| $F$ /(\$/h) | - | - | 8199.11 | 8174.34 | 8175.99 | 8145.06 | 8144.15 | 8150.50 | 8143.41 |
| $IMP_F$ /% | - | - | 0 | 0.302 | 0.282 | 0.659 | 0.670 | 0.593 | 0.679 |
| $V_{de}$ /p.u. | - | - | 0.0077 | 0.0073 | 0.0070 | 0.0056 | 0.0056 | 0.0055 | 0.0057 |
| $IMP_V$ /% | - | - | 0 | 5.752 | 8.999 | 27.358 | 27.074 | 28.365 | 26.339 |

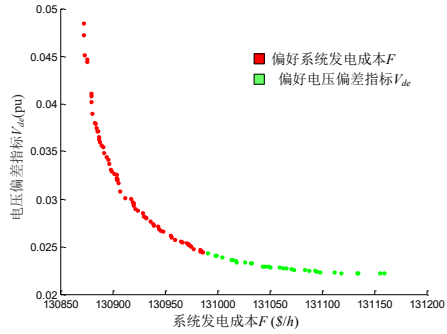

图 A8　IEEE 118 节点算例 Pareto 最优集聚类后的分布
Fig.A8　Distribution of optimal solutions of IEEE 118-bus

表 A10　IEEE 118 节点各工况优化结果对比
Table A10　Comparison of results for 118-bus system

| 变量 | 最大值 | 最小值 | Case 0 | Case 1 | Case 2 | Case 3 |
|---|---|---|---|---|---|---|
| $P_{s,1}$ /p.u. | 1.00 | -1.00 | - | 0.244 | 0.587 | 0.562 |
| $Q_{s,1}$ /p.u. | 1.00 | -1.00 | - | 0.104 | -0.013 | -0.001 |
| $U_{dc,1}$ /p.u. | 1.06 | 0.94 | - | 1.003 | 1.001 | 1.025 |
| $P_{s,2}$ /p.u. | 1.00 | -1.00 | - | -0.234 | -0.571 | -0.566 |
| $Q_{s,2}$ /p.u. | 1.00 | -1.00 | - | 0.068 | 0.085 | 0.089 |
| $U_{dc,2}$ /p.u. | 1.06 | 0.94 | - | 1.000 | 0.998 | 1.025 |
| $P_{s,3}$ /p.u. | 1.00 | -1.00 | - | - | -0.015 | 0.006 |
| $Q_{s,3}$ /p.u. | 1.00 | -1.00 | - | - | 0.000 | 0.004 |
| $U_{dc,3}$ /p.u. | 1.06 | 0.94 | - | - | 1.000 | 1.023 |
| $F$ /(\$/h) | - | - | 131138 | 131110 | 130939 | 130958 |
| $IMP_F$ /% | | | 0 | 0.021 | 0.152 | 0.137 |
| $V_{de}$ /p.u. | - | - | 0.0275 | 0.0273 | 0.0272 | 0.0272 |
| $IMP_V$ /% | | | 0 | 0.791 | 1.043 | 1.192 |